\input amstex
\input amsppt.sty
\magnification=\magstep1
\hsize=33truecc
\vsize=22.2truecm
\baselineskip=16truept
\NoBlackBoxes
\TagsOnRight \pageno=1 \nologo
\def\Z{\Bbb Z}
\def\N{\Bbb N}

\def\l{\left}
\def\r{\right}
\def\bg{\bigg}
\def\({\bg(}
\def\[{\bg\lfloor}
\def\){\bg)}
\def\]{\bg\rfloor}
\def\t{\text}
\def\f{\frac}

\def\p{\ (\roman{mod}\ p)}

\def\sm{\setminus}

\def\bi{\binom}
\def\eq{\equiv}

\def\ls{\leqslant}
\def\gs{\geqslant}
\def\mo{\roman{mod}}

\def\al{\alpha}
\def\da{\delta}

\def\Remark{\medskip\noindent{\it  Remark}}

\def\Ack{\medskip\noindent {\bf Acknowledgment}}
\hbox {Finite Fields Appl. 22(2013), 24--44.}\bigskip
\topmatter
\title Supercongruences involving products of two binomial coefficients\endtitle
\author Zhi-Wei Sun\endauthor
\leftheadtext{Zhi-Wei Sun} \rightheadtext{Supercongruences involving products of binomial coefficients}
\affil Department of Mathematics, Nanjing University\\
 Nanjing 210093, People's Republic of China
  \\ E-mail: zwsun\@nju.edu.cn
  \\ {\tt http://math.nju.edu.cn/$\sim$zwsun}
\endaffil
\abstract In this paper we deduce some new supercongruences
modulo powers of a prime $p>3$. Let $d\in\{0,1,\ldots,(p-1)/2\}$. We show that
 $$\sum_{k=0}^{(p-1)/2}\f{\bi{2k}k\bi{2k}{k+d}}{8^k}\eq0\ (\mo\ p)\ \ \ \t{if}\ d\eq \f{p+1}2\ (\mo\ 2),$$
and $$\sum_{k=0}^{(p-1)/2}\f{\bi{2k}k\bi{2k}{k+d}}{16^k}
\eq\l(\f{-1}p\r)+p^2\f{(-1)^d}4E_{p-3}\l(d+\f12\r)\pmod{p^3},$$
where $E_{p-3}(x)$ denotes the Euler polynomial of degree $p-3$, and $(-)$ stands for the Legendre symbol.
The paper also contains some other results such as
$$\sum_{k=0}^{p-1}k^{(1+(\f{-1}p))/2}\f{\bi{6k}{3k}\bi{3k}k}{864^k}\eq0\pmod{p^2}.$$
\endabstract
\thanks 2010 {\it Mathematics Subject Classification}. Primary 11B65;
Secondary 05A10, 11A07, 11B39, 11B68, 11E25, 12E20.
\newline\indent Supported by the National Natural Science
Foundation (grant 11171140) of China and the PAPD of Jiangsu Higher
Education Institutions.
\endthanks
\endtopmatter
\document

\heading{1. Introduction}\endheading

Let $p$ be an odd prime and let $(\f{\cdot}p)$ be the Legendre symbol. For each $d\in\N=\{0,1,\ldots\}$
and any rational $p$-adic integer $\lambda$, we define
$$a_p^{(d)}(\lambda):=\sum_{x=0}^{p-1}x^d\l(\f{x(x-1)(x-\lambda)}p\r).\tag1.1$$
Note that $a_p^{(0)}(\lambda)$ arises naturally from counting the number of points on the cubic curve
${\Bbb E}_p(\lambda):\ y^2=x(x-1)(x-\bar \lambda)$
over the finite field $\Bbb F_p=\Z/p\Z$, where $\bar\lambda$ is the residue class $\lambda\ (\mo\ p)$.

The following theorem in the case $d=0$ is a known result (cf. S. Ahlgren [A, Theorem 2]).

\proclaim{Theorem 1.1} Let $p$ be an odd prime and let $d\in\{0,\ldots,(p-1)/2\}$. Then, for any rational
$p$-adic integer $\lambda$ we have
$$a_p^{(d)}(\lambda)\eq (-1)^{(p+1)/2}\f{\lambda^d}{4^d}\sum_{k=0}^{(p-1)/2}\f{\bi{2k}k\bi{2(k+d)}{k+d}}{16^k}\lambda^k
-\da_{d,(p-1)/2}\pmod{p},\tag1.2$$
where the Kronecker symbol $\da_{s,t}$ takes $1$ or $0$ according as $s=t$ or not.
\endproclaim

\Remark\ 1.1. Let  $d\in\{0,\ldots,(p-1)/2\}$ with $p$ an odd prime. Clearly
$$a_p^{(d)}(1)=\sum_{x=0}^{p-1}x^d\l(\f xp\r)-1\eq\sum_{x=1}^{p-1}x^{d+(p-1)/2}-1\eq-\da_{d,(p-1)/2}-1\ (\mo\ p).$$
Thus (1.2) with $\lambda=1$ gives the congruence
$$\sum_{k=0}^{(p-1)/2}\f{\bi{2k}k\bi{2k+2d}{k+d}}{16^k}\eq 4^d\l(\f{-1}p\r)\pmod{p}.$$
Soon we will see that this congruence even holds modulo $p^2$.
\medskip

Recall that the Euler numbers $E_0,E_1,E_2,\ldots$ are integers defined
by $$E_0=1\ \ \t{and}\ \ \sum^n\Sb k=0\\2\mid k\endSb \bi nk
E_{n-k}=0\ \ \ \t{for}\ n\in\Z^+=\{1,2,3,\ldots\}.$$
For each $n\in\N$, the Euler polynomial of degree $n$ is given by
$$ E_n(x)=\sum_{k=0}^n\bi{n}{k}\frac{E_k}{2^k}\left(x-\frac{1}{2}\right)^{n-k}.$$
Clearly $E_n(1/2)=E_n/2^n$.

Now we state our second theorem.

 \proclaim{Theorem 1.2} Let $p>3$ be a prime and let
$d\in\{0,\ldots,(p-1)/2\}$. Then
$$\sum_{k=0}^{(p-1)/2}\f{\bi{2k}k\bi{2k}{k+d}}{16^k}\eq\l(\f{-1}p\r)+p^2\f{(-1)^d}4E_{p-3}\l(d+\f12\r)\pmod{p^3}.\tag1.3$$
\endproclaim

 \Remark\ 1.2. Let $p>3$ be a prime. The supercongruence
$$\sum_{k=0}^{(p-1)/2}\f{\bi{2k}k^2}{16^k}\eq\l(\f{-1}p\r)\pmod{p^2}$$
  was a conjecture of Rodriguez-Villegas [RV] confirmed by E. Mortenson [Mo1]
 via an advanced tool involving the $p$-adic Gamma function
and the Gross-Koblitz formula for character sums.  (1.3) with $d=0$ yields the congruence
$$\sum_{k=0}^{(p-1)/2}\f{\bi{2k}k^2}{16^k}\eq\l(\f{-1}p\r)+p^2E_{p-3}\pmod{p^3}$$
which was first proved in [S4] with the help of the software {\tt Sigma}.

\proclaim{Corollary 1.1} Let $p>3$ be a prime. For any $d=0,\ldots,(p-1)/2$, we have
 $$\sum_{k=0}^{(p-1)/2}\f{\bi{2k}k\bi{2k+2d}{k+d}}{16^k}\eq 4^d\l(\f{-1}p\r)\pmod{p^2}.\tag1.4$$
\endproclaim

\medskip
Let $p\eq1\ (\mo\ 4)$ be a prime. It is well known that $p=x^2+y^2$ for some $x,y\in\Z$ with $x\eq1\ (\mo\ 4)$.
A celebrated result of Gauss asserts that
$\bi{(p-1)/2}{(p-1)/4}\eq2x\ (\mo\ p)$
(see, e.g., [BEW, (9.0.1)]).
This was refined in [CDE] as follows:
$$\bi{(p-1)/2}{(p-1)/4}\eq\f{2^{p-1}+1}2\l(2x-\f p{2x}\r)\pmod{p^2}.$$
Recently, J. B. Cosgrave and K. Dilcher [CD] even determined $\bi{(p-1)/2}{(p-1)/4}$ mod $p^3$.
Recall that $p=x^2+y^2$ with $x\eq1\ (\mo\ 4)$.
Z.-H. Sun [Su] confirmed the author's following conjecture (cf. [S3, Conjecture 5.5]):
$$\align\sum_{k=0}^{(p-1)/2}\f{\bi{2k}k^2}{8^k}\eq&\sum_{k=0}^{(p-1)/2}\f{\bi{2k}k^2}{(-16)^k}
\eq\l(\f 2p\r)\sum_{k=0}^{(p-1)/2}\f{\bi{2k}k^2}{32^k}
\\\eq&\l(\f 2p\r)\l(2x-\f p{2x}\r)\ (\mo\ p^2).
\endalign$$
In [S5] the author showed that
$$\align&\sum_{k=0}^{(p-1)/2}\f{\bi{2k}{k}^2}{(k+1)8^k}\eq-2\sum_{k=0}^{p-1}\f{k\bi{2k}k^2}{8^k}
\eq\f12\sum_{k=0}^{(p-1)/2}\f{\bi{2k}{k}^2}{(k+1)(-16)^k}
\\\eq&-4\sum_{k=0}^{(p-1)/2}\f{k\bi{2k}k^2}{(-16)^k}\eq\l(\f 2p\r)\l(2x-\f px\r)\pmod{p^2}.
\endalign$$
Note that those integers $\bi{2k}k/(k+1)=\bi{2k}k-\bi{2k}{k+1}$ are called Catalan numbers and they
 occur naturally in many enumeration problems in combinatorics (see, e.g., [St, pp.\,219--229]).

Motivated by (1.2) in the cases $\lambda=-1,2$ we obtain the following result.
\proclaim{Theorem 1.3} {\rm (i)} If $p\eq3\pmod 4$ is a prime, then
$$\aligned&\sum_{k=0}^{(p-1)/2}\f{\bi{2k}{k}^2}{(k+1)8^k}\eq-2\sum_{k=0}^{(p-1)/2}\f{k\bi{2k}k^2}{8^k}
\\\eq&-\f12\sum_{k=0}^{(p-1)/2}\f{\bi{2k}{k}^2}{(k+1)(-16)^k}\eq4\sum_{k=0}^{(p-1)/2}\f{k\bi{2k}k^2}{(-16)^k}
\\\eq&\f{(-1)^{(p+1)/4}}2\bi{(p+1)/2}{(p+1)/4}\pmod{p}
\endaligned\tag1.5$$
and
$$\sum_{k=0}^{(p-1)/2}\f{\bi{2k}k^2}{8^k}\eq-\sum_{k=0}^{(p-1)/2}\f{\bi{2k}k^2}{(-16)^k}
\eq\f{(-1)^{(p+1)/4}\, 2p}{\bi{(p+1)/2}{(p+1)/4}}\pmod{p^2}.\tag1.6$$

{\rm (ii)} For any odd prime $p$, we have
$$\sum_{k=0}^{(p-1)/2}\f{\bi{2k}k\bi{2k}{k+d}}{8^k}\eq0\pmod{p}\tag1.7$$
for all $d\in\{0,\ldots,(p-1)/2\}$ with $d\eq (p+1)/2\ (\mo\ 2)$.
\endproclaim
\Remark\ 1.3. In 2009 the author conjectured that $\sum_{k=0}^{p-1}\bi{2k}k\bi{2k}{k+1}/8^k\eq0\ (\mo\ p)$
for any prime $p\eq1\ (\mo\ 4)$ and this was confirmed by his student Yong Zhang in her PhD thesis.
\medskip

Besides (1.4) with $d=0$, Rodriguez-Villegas [RV] also made the
following similar conjectures (confirmed in [Mo2]) on supercongruences with $p$ a prime greater than 3:
$$\align \sum_{k=0}^{p-1}\f{\bi{3k}k\bi{2k}k}{27^k}\eq&\l(\f{p}3\r)\ (\mo\ p^2),\tag1.8
\\\sum_{k=0}^{p-1}\f{\bi{4k}{2k}\bi{2k}k}{64^k}\eq&\l(\f{-2}p\r)\ (\mo\ p^2),\tag1.9
\\\sum_{k=0}^{p-1}\f{\bi{6k}{3k}\bi{3k}k}{432^k}\eq&\l(\f{-1}p\r)\ (\mo\ p^2).\tag1.10
\endalign$$
Note that the denominators 27, 64, 432 come from the following observation via the Stirling formula:
$$\gather\binom{3k}k\binom{2k}k\sim\f{\sqrt3\times 27^k}{2k\pi},\ \
\binom{4k}{2k}\binom{2k}{k}\sim\f{64^k}{\sqrt2k\pi},\ \ \binom{6k}{3k}\binom{3k}{k}\sim\f{432^k}{2k\pi}.
\endgather$$
Up to now no simple proofs of (1.8)-(1.10) have been found.

Motivated by the work in [PS] and [ST], the author [S2] determined $\sum_{k=0}^{p-1}\bi{2k}k/m^k$ modulo $p^2$
in terms of Lucas sequences, where
$p$ is an odd prime and $m$ is any integer not divisible by $p$.
 In [S3] and [S4] the author posed many conjectures on sums of terms involving central binomial coefficients.

 For a sequence of $(a_n)_{n\in\N}$ of numbers, as in [S1] we introduce its dual sequence $(a_n^*)_{n\in\N}$
by defining
$$a_n^*:=\sum_{k=0}^n\bi nk(-1)^ka_k\ \ \  (n=0,1,2,\ldots).$$
It is well-known that $(a_n^*)^*=a_n$ for all $n\in\N$ (see, e.g., (5.48) of [GKP, p.\,192]).
For Bernoulli numbers $B_0,B_1,B_2,\ldots$, the sequence $((-1)^nB_n)_{n\in\N}$ is self-dual.

\proclaim{Theorem 1.4} Let $p>3$ be a prime and let $(a_n)_{n\in\N}$ be any sequence of $p$-adic integers. Then we have
$$\align \sum_{k=0}^{p-1}\f{\bi{3k}k\bi{2k}k}{27^k}a_k\eq&\l(\f p3\r)\sum_{k=0}^{p-1}\f{\bi{3k}k\bi{2k}k}{27^k}a_k^*\pmod{p^2},\tag1.11
\\\sum_{k=0}^{p-1}\f{\bi{4k}{2k}\bi{2k}k}{64^k}a_k\eq&\l(\f {-2}p\r)\sum_{k=0}^{p-1}\f{\bi{4k}{2k}\bi{2k}k}{64^k}a_k^*\pmod{p^2},\tag1.12
\\\sum_{k=0}^{p-1}\f{\bi{6k}{3k}\bi{3k}k}{432^k}a_k\eq&\l(\f {-1}p\r)\sum_{k=0}^{p-1}\f{\bi{6k}{3k}\bi{3k}k}{432^k}a_k^*\pmod{p^2}.\tag1.13
\endalign$$
\endproclaim
\Remark\ 1.4. Z.-H. Sun [Su] recently proved that
$$\sum_{k=0}^{(p-1)/2}\f{\bi{2k}k^2}{16^k}\l(a_k-\l(\f {-1}p\r)a_k^*\r)\eq0\ (\mo\ p^2)$$
for any odd prime $p$ via Legendre polynomials.
We can also show, for any prime $p>3$, the following result similar to (1.3) and (1.4):
If $d\in\{0,\ldots,\lfloor p/3\rfloor\}$ then
$$\f1{4^d}\sum_{k=0}^{(p-1)/2}\f{\bi{3k}k\bi{2k+2d}{k+d}}{27^k}
\eq\sum_{k=0}^{(p-1)/2}\f{\bi{3k}k\bi{2k}{k+d}}{27^k}\eq\l(\f p3\r)\pmod{p};$$
if $d\in\{0,\ldots,\lfloor p/4\rfloor\}$ then
$$\f1{4^d}\sum_{k=0}^{(p-1)/2}\f{\bi{4k}{2k}\bi{2k+2d}{k+d}}{64^k}
\eq\sum_{k=0}^{(p-1)/2}\f{\bi{4k}{2k}\bi{2k}{k+d}}{64^k}\eq\l(\f {-2}p\r)\pmod{p}.$$

\medskip

Let $p$ be a prime and let $f(x)\in\Bbb F_p[x]$ with $\deg(f)<p$. Then $f(x)$ is identically zero if $f(a)=0$ for all $a\in \Bbb F_p$.
Thus Theorem 1.4 has the following consequence since $(1-x)^k=\sum_{j=0}^k\bi kj(-1)^jx^j$ for any $k\in\N$.

\proclaim{Corollary 1.2} Let $p>3$ be a prime and let $\Z_p$ be the ring of $p$-adic integers. Then, in the ring $\Z_p[x]$ we have
$$\align \sum_{k=0}^{p-1}\f{\bi{3k}k\bi{2k}k}{27^k}\l(x^k-\l(\f p3\r)(1-x)^k\r)\eq&0\pmod{p^2},\tag1.14
\\\sum_{k=0}^{p-1}\f{\bi{4k}{2k}\bi{2k}k}{64^k}\l(x^k-\l(\f {-2}p\r)(1-x)^k\r)\eq&0\pmod{p^2},\tag1.15
\\\sum_{k=0}^{p-1}\f{\bi{6k}{3k}\bi{3k}k}{432^k}\l(x^k-\l(\f {-1}p\r)(1-x)^k\r)\eq&0\pmod{p^2}.\tag1.16
\endalign$$
Also,
$$\align \sum_{k=1}^{p-1}\f{k\bi{3k}k\bi{2k}k}{27^k}\l(x^{k-1}+\l(\f p3\r)(1-x)^{k-1}\r)\eq&0\pmod{p^2},\tag1.17
\\\sum_{k=1}^{p-1}\f{k\bi{4k}{2k}\bi{2k}k}{64^k}\l(x^{k-1}+\l(\f {-2}p\r)(1-x)^{k-1}\r)\eq&0\pmod{p^2},\tag1.18
\\\sum_{k=1}^{p-1}\f{k\bi{6k}{3k}\bi{3k}k}{432^k}\l(x^{k-1}+\l(\f {-1}p\r)(1-x)^{k-1}\r)\eq&0\pmod{p^2}.\tag1.19
\endalign$$
\endproclaim
\Remark\ 1.5. (1.17)-(1.19) can be easily deduced from (1.14)-(1.16) by taking derivatives.
Z.-H. Sun [Su, Theorem 2.4] noted that for any prime $p>3$ we have
$$\sum_{k=0}^{\lfloor p/3\rfloor}\f{\bi{3k}k\bi{2k}k}{27^k}(x^k-(-1)^{\lfloor p/3\rfloor}(1-x)^k)\eq0\pmod{p}.$$

\medskip

Taking $x=1/2$ in (1.14)-(1.19) we immediately get the following result.

\proclaim{Corollary 1.3} Let $p>3$ be a prime. Then
$$\align\sum_{k=0}^{p-1}\f{k\bi{3k}k\bi{2k}k}{54^k}\eq&0\ (\mo\ p^2)\ \ \t{if}\ p\eq1\ (\mo\ 3),
\\\sum_{k=0}^{p-1}\f{\bi{3k}k\bi{2k}k}{54^k}\eq&0\ (\mo\ p^2)\ \ \t{if}\ p\eq2\ (\mo\ 3);
\\\sum_{k=0}^{p-1}\f{k\bi{4k}{2k}\bi{2k}k}{128^k}\eq&0\ (\mo\ p^2)\ \ \t{if}\ p\eq1,3\ (\mo\ 8),
\\\sum_{k=0}^{p-1}\f{\bi{4k}{2k}\bi{2k}k}{128^k}\eq&0\ (\mo\ p^2)\ \ \t{if}\ p\eq5,7\ (\mo\ 8);
\\\sum_{k=0}^{p-1}\f{k\bi{6k}{3k}\bi{3k}k}{864^k}\eq&0\ (\mo\ p^2)\ \ \t{if}\ p\eq1\ (\mo\ 4),
\\\sum_{k=0}^{p-1}\f{\bi{6k}{3k}\bi{3k}k}{864^k}\eq&0\ (\mo\ p^2)\ \ \t{if}\ p\eq3\ (\mo\ 4).
\endalign$$
\endproclaim
\Remark\ 1.6. The first and the second congruences mod $p$ were obtained by Z.-H. Sun [Su].
The author [S4] and Z.-H. Sun [Su] conjectured the first and the second congruences respectively.
Inputting
\newline {\tt FullSimplify[Sum[k*Binomial[3k,k]Binomial[2k,k]/$54^\wedge$k,$\{$k,0,Infty$\}$]]},
\newline we obtain from {\tt Mathematica} the exact result
$$\sum_{k=0}^\infty\f{k\bi{3k}k\bi{2k}k}{54^k}=\f{\sqrt{\pi}}{9\Gamma(\f43)\Gamma(\f 76)},$$
which should follow from certain algorithm hidden in {\tt Mathematica}.

\medskip
(1.14) and (1.17) in the case $x=9/8$, and (1.15) and (1.18) in the cases $x=4/3,\,8/9,\,64/63$, yield the following result.

\proclaim{Corollary 1.4} Let $p>3$ be a prime. Then
$$\align \sum_{k=0}^{p-1}\f{\bi{3k}k\bi{2k}k}{24^k}\eq&\l(\f p3\r)\sum_{k=0}^{p-1}\f{\bi{3k}k\bi{2k}k}{(-216)^k}\pmod{p^2},\tag1.20
\\ \sum_{k=0}^{p-1}\f{k\bi{3k}k\bi{2k}k}{24^k}\eq&9\l(\f p3\r)\sum_{k=0}^{p-1}\f{k\bi{3k}k\bi{2k}k}{(-216)^k}\pmod{p^2}.\tag1.21
\endalign$$
Also,
$$\align\sum_{k=0}^{p-1}\f{\bi{4k}{2k}\bi{2k}k}{48^k}\eq&\l(\f {-2}p\r)\sum_{k=0}^{p-1}\f{\bi{4k}{2k}\bi{2k}k}{(-192)^k}\ (\mo\ p^2),
\\ \sum_{k=0}^{p-1}\f{k\bi{4k}{2k}\bi{2k}k}{48^k}\eq&4\l(\f {-2}p\r)\sum_{k=0}^{p-1}\f{k\bi{3k}k\bi{2k}k}{(-192)^k}\ (\mo\ p^2);
\\\sum_{k=0}^{p-1}\f{\bi{4k}{2k}\bi{2k}k}{72^k}\eq&\l(\f {-2}p\r)\sum_{k=0}^{p-1}\f{\bi{4k}{2k}\bi{2k}k}{576^k}\ (\mo\ p^2),
\\ \sum_{k=0}^{p-1}\f{k\bi{4k}{2k}\bi{2k}k}{72^k}\eq&-8\l(\f {-2}p\r)\sum_{k=0}^{p-1}\f{k\bi{3k}k\bi{2k}k}{576^k}\ (\mo\ p^2).
\endalign$$
If $p\not=7$, then
$$\align\sum_{k=0}^{p-1}\f{\bi{4k}{2k}\bi{2k}k}{63^k}\eq&\l(\f {-2}p\r)\sum_{k=0}^{p-1}\f{\bi{4k}{2k}\bi{2k}k}{(-4032)^k}\pmod{p^2},
\\ \sum_{k=0}^{p-1}\f{k\bi{4k}{2k}\bi{2k}k}{63^k}\eq&64\l(\f {-2}p\r)\sum_{k=0}^{p-1}\f{k\bi{3k}k\bi{2k}k}{(-4032)^k}\pmod{p^2}.
\endalign$$
\endproclaim
\Remark\ 1.7. Let $p>3$ be a prime. In [S4, Conjecture 5.13] the author conjectured that
$$\sum_{k=0}^{p-1}\f{\bi{2k}k\bi{3k}k}{24^k}\eq\l(\f p3\r)\sum_{k=0}^{p-1}\f{\bi{2k}k\bi{3k}k}{(-216)^k}
\eq\cases\bi{2(p-1)/3}{(p-1)/3}\ (\mo\ p^2)&\t{if}\ p\eq1\ (\mo\ 3),
\\p/\bi{2(p+1)/3}{(p+1)/3}\ (\mo\ p^2)&\t{if}\ p\eq 2\ (\mo\ 3).\endcases$$
The author [S4] also made conjectures on
$\sum_{k=0}^{p-1}\bi{4k}{2k}\bi{2k}k/m^k$ modulo $p^2$
with $m=48,63,72,128$.
\medskip

For any prime $p>3$ and integer $m\not\eq0\pmod p$, we have
$$\align\sum_{k=0}^{p-1}\f{\bi{3k}{k}\bi{2k}k}{(k+1)m^k}\eq& p+\f{m-27}{6}\sum_{k=0}^{p-1}\f{k\bi{3k}{k}\bi{2k}k}{m^k}\pmod{p^2},\tag1.22
\\\sum_{k=0}^{p-1}\f{\bi{4k}{2k}\bi{2k}k}{(k+1)m^k}\eq& p+\f{m-64}{12}\sum_{k=0}^{p-1}\f{k\bi{4k}{2k}\bi{2k}k}{m^k}\pmod{p^2},\tag1.23
\\\sum_{k=0}^{p-1}\f{\bi{6k}{3k}\bi{3k}k}{(k+1)m^k}\eq& p+\f{m-432}{60}\sum_{k=0}^{p-1}\f{k\bi{6k}{3k}\bi{3k}k}{m^k}\pmod{p^2},\tag1.24
\endalign$$
due to the identities
$$\align\sum_{k=0}^{n-1}\l(\f{6\bi{2k}k}{k+1}+(27-m)k\bi{2k}k\r)\f{\bi{3k}{k}}{m^k}=&\f{n}{m^{n-1}}\bi{2n}{n}\bi{3n}n,
\\\sum_{k=0}^{n-1}\l(\f{12\bi{2k}k}{k+1}+(64-m)k\bi{2k}k\r)\f{\bi{4k}{2k}}{m^k}=&\f{n}{m^{n-1}}\bi{4n}{2n}\bi{2n}n,
\\\sum_{k=0}^{n-1}\l(\f{60}{k+1}+(432-m)k\r)\f{\bi{6k}{3k}\bi{3k}k}{m^k}=&\f{n}{m^{n-1}}\bi{6n}{3n}\bi{3n}n,
\endalign$$
which can be easily proved by induction on $n$. So, the following result follows from Corollary 1.3 and
(1.21).

\proclaim{Corollary 1.5} Let $p>3$ be a prime. Then
$$\align\sum_{k=0}^{p-1}\f{\bi{3k}k\bi{2k}k}{(k+1)54^k}\eq& p\ \ (\mo\ p^2)\ \ \t{if}\ p\eq1\ (\mo\ 3),\tag1.25
\\\sum_{k=0}^{p-1}\f{\bi{4k}{2k}\bi{2k}k}{(k+1)128^k}\eq& p\ \ (\mo\ p^2)\ \ \t{if}\ p\eq1,3\ (\mo\ 8),\tag1.26
\\\sum_{k=0}^{p-1}\f{\bi{6k}{3k}\bi{3k}k}{(k+1)864^k}\eq& p\ \ (\mo\ p^2)\ \ \t{if}\ p\eq1\ (\mo\ 4).\tag1.27
\endalign$$
We also have
$$\sum_{k=0}^{p-1}\f{\bi{3k}k\bi{2k}k}{(k+1)24^k}
\eq p+\f19\l(\f p3\r)\(\sum_{k=0}^{p-1}\f{\bi{3k}k\bi{2k}k}{(k+1)(-216)^k}-p\)\ (\mo\ p^2).\tag1.28$$
\endproclaim
\Remark\ 1.8. Similar to the identity in Remark 1.6, {\tt Mathematica} (version 7) also yields
 $$\sum_{k=0}^\infty\f{\bi{3k}k\bi{2k}k}{(k+1)54^k}=\f{3\sqrt{\pi}}{\Gamma(\f43)\Gamma(\f16)},
\  \ \sum_{k=0}^\infty\f{\bi{4k}{2k}\bi{2k}k}{(k+1)128^k}=\f{4\sqrt{\pi}}{\Gamma(\f18)\Gamma(\f{11}8)},$$
 and
$$\sum_{k=0}^\infty\f{\bi{6k}{3k}\bi{3k}k}{(k+1)864^k}=\f{6\sqrt{\pi}}{\Gamma(\f1{12})\Gamma(\f{17}{12})}.$$

\proclaim{Theorem 1.5} Let $p>3$ be a prime. Then
$$\aligned&\sum_{k=0}^{p-1}\f{k\bi{4k}{2k}\bi{2k}k}{72^k}\eq\f32\sum_{k=0}^{p-1}\f{\bi{4k}{2k}\bi{2k}k}{(k+1)72^k}
\\\eq &\cases(\f 6p)x\ (\mo\ p)&\t{if}\ p=x^2+y^2\ \t{with}\ x\eq1\ (\mo\ 4),
\\\f 34(\f 6p)\bi{(p+1)/2}{(p+1)/4}\ (\mo\ p)&\t{if}\ p\eq3\pmod4.
\endcases\endaligned\tag1.29$$
\endproclaim

Let $A$ and $B$ be integers. The Lucas sequences $u_n=u_n(A,B)\ (n\in\N)$ and $v_n=v_n(A,B)\ (n\in\N)$ are defined as follows:
$$\align &u_0=0,\ u_1=1,\ \t{and}\ u_{n+1}=Au_n-Bu_{n-1}\ \t{for}\ n\in\Z^+;
\\&v_0=2,\ v_1=A,\ \t{and}\ v_{n+1}=Av_n-Bv_{n-1}\ \t{for}\ n\in\Z^+.
\endalign$$
When $\Delta=A^2-4B=0$, by induction we see that $u_n(A,B)=n(A/2)^{n-1}$ and $v_n(A,B)=2(A/2)^n$ for all $n\in\Z^+$.
Our following theorem is an analogue of Corollary 1.3 involving Lucas sequences with $\Delta\not=0$.

\proclaim{Theorem 1.6} Let $A,B\in\Z$ with $A\not=0$ and $A^2\not=4B$, and let $u_k=u_k(A,B)$ and $v_k=v_k(A,B)$ for all $k\in\N$.
Let $p>3$ be a prime with $p\nmid A$.

{\rm (i)} If $p\eq1\ (\mo\ 3)$, then
$$\sum_{k=0}^{p-1}\f{\bi{3k}k\bi{2k}k}{(27A)^k}u_k\eq \sum_{k=1}^{p-1}\f{k\bi{3k}k\bi{2k}k}{(27A)^k}v_{k-1}\eq0\pmod{p^2}.\tag1.30$$
If $p\eq2\ (\mo\ 3)$, then
$$\sum_{k=0}^{p-1}\f{\bi{3k}k\bi{2k}k}{(27A)^k}v_k\eq \sum_{k=1}^{p-1}\f{k\bi{3k}k\bi{2k}k}{(27A)^k}u_{k-1}\eq0\pmod{p^2}.\tag1.31$$

{\rm (ii)} If $p\eq1,3\ (\mo\ 8)$, then
$$\sum_{k=0}^{p-1}\f{\bi{4k}{2k}\bi{2k}k}{(64A)^k}u_k\eq \sum_{k=1}^{p-1}\f{k\bi{4k}{2k}\bi{2k}k}{(64A)^k}v_{k-1}\eq0\pmod{p^2}.\tag1.32$$
If $p\eq5,7\ (\mo\ 8)$, then
$$\sum_{k=0}^{p-1}\f{\bi{4k}{2k}\bi{2k}k}{(64A)^k}v_k\eq \sum_{k=1}^{p-1}\f{k\bi{4k}{2k}\bi{2k}k}{(64A)^k}u_{k-1}\eq0\pmod{p^2}.\tag1.33$$

{\rm (iii)} If $p\eq1\ (\mo\ 4)$, then
$$\sum_{k=0}^{p-1}\f{\bi{6k}{3k}\bi{3k}k}{(432A)^k}u_k\eq \sum_{k=1}^{p-1}\f{k\bi{6k}{3k}\bi{3k}k}{(432A)^k}v_{k-1}\eq0\pmod{p^2}.\tag1.34$$
If $p\eq3\ (\mo\ 4)$, then
$$\sum_{k=0}^{p-1}\f{\bi{6k}{3k}\bi{2k}k}{(432A)^k}v_k\eq \sum_{k=1}^{p-1}\f{k\bi{6k}{3k}\bi{2k}k}{(432A)^k}u_{k-1}\eq0\pmod{p^2}.\tag1.35$$
\endproclaim

We will not list corollaries of Theorem 1.6 with respect to some special Lucas sequences like the Fibonacci sequence $F_n=u_n(1,-1)\ (n\in\N)$
and its companion $L_n=v_n(1,-1)\ (n\in\N)$.

In the next section we are going to show Theorems 1.1-1.3 and Corollary 1.1. Sections 3 and 4 are devoted to our proofs of
Theorem 1.4 and Theorems 1.5-1.6 respectively.

\heading{2. Proofs of Theorems 1.1-1.3 and Corollary 1.1}\endheading

\medskip
\noindent{\it Proof of Theorem 1.1}. Set $n=(p-1)/2$. Then
$$\align a_p^{(d)}(\lambda)\eq&\sum_{k=0}^{p-1}x^d\l(x(x-1)(x-\lambda)\r)^{n}
\\=&\sum_{k=0}^{p-1}x^{n+d}\sum_{k=0}^{n}\bi{n}k(-1)^{n-k}x^k\sum_{l=0}^{n}\bi{n}l(-\lambda)^lx^{n-l}
\\=&\sum_{k,l=0}^{n}\bi{n}k\bi{n}l(-1)^{n-k}(-\lambda)^l\sum_{x=1}^{p-1}x^{p-1+d+k-l}
\\\eq&\sum_{k=0}^n\bi nk(-1)^{n-k}\sum\Sb 0\ls l\ls n\\p-1\mid l-(d+k)\endSb\bi nl(-\lambda)^l(p-1)
\\\eq&-\sum_{k=0}^n\bi nk(-1)^{n-k}\bi n{d+k}(-\lambda)^{d+k}-\da_{d,n}\bi n0(-\lambda)^0\ (\mo\ p).
\endalign$$
Since
$$\bi {(p-1)/2}k\eq\bi{-1/2}k=\f{\bi{2k}k}{(-4)^k}\pmod{p}\quad\t{for all}\ k=0,\ldots,p-1,\tag2.1$$
we immediately obtain (1.2) from the above. \qed

\medskip
\noindent{\it Proof of Theorem 1.2}. By induction, we have
$$\sum_{k=0}^n\f{\bi{2k}k}{16^k}\l(\bi{2k}{k+m}-\bi{2k}{k+m+1}\r)=\f{2n+1}{(2m+1)16^n}\bi{2n}n\bi{2n+1}{n-m}\tag2.2$$
for each $n=m,m+1,\ldots$, where $m\in\N$.

Set $n=(p-1)/2$. If $0\ls m<n$,  then for the right-hand side $R_m$ of (2.2) we have
$$\align R_m=&\f{p^2}{(2m+1)((p-1)/2-m)4^{p-1}}\bi{p-1}{n}\bi{p-1}{n-m-1}
\\\eq&2p^2\f{(-1)^m}{(2m+1)^2}\ \ \pmod{p^3}
\endalign$$
since
$$\bi{p-1}k=\prod_{0<j\ls k}\f{p-j}j\eq(-1)^k\pmod{p}\quad\t{for all}\ k=0,\ldots,p-1.\tag2.3$$
As $d\ls n$, we have
$$\align &\sum_{k=0}^n\f{\bi{2k}k}{16^k}\l(\bi{2k}{k}-\bi{2k}{k+d}\r)
\\=&\sum_{0\ls m<d}\sum_{k=0}^n\f{\bi{2k}k}{16^k}\l(\bi{2k}{k+m}-\bi{2k}{k+m+1}\r)
\\\eq&2p^2\sum_{0\ls m<d}\f{(-1)^m}{(2m+1)^2}\eq\f{p^2}2\sum_{0\ls m<d}(-1)^m\l(m+\f12\r)^{p-3}
\\=&\f{p^2}4\sum_{0\ls m<d}(-1)^m\(E_{p-3}\l(m+\f12\r)+E_{p-3}\l(m+1+\f12\r)\)
\\=&\f{p^2}4\l(E_{p-3}\l(\f12\r)-(-1)^dE_{p-3}\l(d+\f12\r)\r)\pmod{p^3}.
\endalign$$
Note that
$$\sum_{k=0}^n\f{\bi{2k}k\bi{2k}{k+n}}{16^k}=\f{\bi{2n}n}{16^n}=\f{\bi{p-1}{(p-1)/2}}{4^{p-1}}\eq\l(\f{-1}p\r)=(-1)^n\pmod{p^3}$$
by Morley's congruence ([M]), and that
$$E_{p-3}\l(n+\f12\r)=E_{p-3}\l(\f p2\r)\eq E_{p-3}(0)=0\pmod{p}.$$
(It is well known that $E_{2k}(0)=0$ for all $k\in\Z^+$.)
Therefore
$$\sum_{k=0}^n\f{\bi{2k}k^2}{16^k}-(-1)^n\eq\f{p^2}4E_{p-3}\l(\f12\r)\eq p^2E_{p-3}\pmod{p^3}$$
and hence (1.3) follows from the above. \qed

\medskip
\noindent{\it Proof of Corollary 1.1}. For $k=0,1,\ldots$, we have
 $$\bi{2k+2d}{k+d}=\sum_{c=-d}^d\bi{2k}{k+c}\bi{2d}{d-c}$$
 by the Chu-Vandermonde identity (cf. [GKP, p.\,169]).
(Note that $\bi{2k}{k+c}$ is regarded as zero if $k+c<0$.)
 In view of this and (1.3),
 $$\align\sum_{k=0}^{(p-1)/2}\f{\bi{2k}k\bi{2k+2d}{k+d}}{16^k}=&\sum_{c=-d}^d\bi{2d}{d-c}\sum_{k=0}^{(p-1)/2}\f{\bi{2k}k\bi{2k}{k+|c|}}{16^k}
\\\eq&\sum_{c=-d}^d\bi{2d}{d-c}\l(\f{-1}p\r)=2^{2d}\l(\f{-1}p\r)\ (\mo\ p^2).
\endalign$$
So (1.4) is valid and we are done. \qed

\medskip
\noindent{\it Proof of Theorem 1.3}. (i) For $m\in\Z\sm\{0\}$ and $n\in\N$ we have the combinatorial identity
$$\sum_{k=0}^n\l(\f{16-m}4k+\f1{k+1}\r)\f{\bi{2k}k^2}{m^k}=\f{(2n+1)^2}{(n+1)m^n}\bi{2n}n^2\tag2.4$$
which can be easily proved by induction on $n$.

Now let $p=2n+1$ be a prime with $p\eq3\ (\mo\ 4)$. Setting $n=(p-1)/2$ we obtain from (2.4) that
$$\sum_{k=0}^n\f{\bi{2k}{k}^2}{(k+1)m^k}\eq\f{m-16}4\sum_{k=0}^n\f{k\bi{2k}k^2}{m^k}\pmod{p^2}\tag2.5$$
for any integer $m\not\eq0\pmod p$.

As $n=(p-1)/2$ is odd, by a result of Z.-H. Sun [Su],
$$\sum_{k=0}^n\f{\bi{2k}k^2}{16^k}(x^k+(1-x)^k)=p^2f(x)$$
for some polynomial $f(x)$ of degree at most $(p-1)/2$ with rational $p$-adic integer coefficients.
In particular,
$$\sum_{k=0}^n\f{\bi{2k}k^2}{8^k}\eq-\sum_{k=0}^n\f{\bi{2k}k^2}{(-16)^k}\pmod{p^2}.\tag2.6$$
By integration,
$$\sum_{k=0}^n\f{\bi{2k}k^2}{(k+1)16^k}x^{k+1}-\sum_{k=0}^n\f{\bi{2k}k^2}{(k+1)16^k}\l((1-x)^{k+1}-1\r)
=p^2\int_0^xf(t)dt.$$
Putting $x=-1$ we obtain
$$-\sum_{k=0}^n\f{\bi{2k}{k}^2}{(k+1)(-16)^k}-\sum_{k=0}^n\f{\bi{2k}{k}^2}{(k+1)16^k}\l(2^{k+1}-1\r)\eq0\pmod{p^2}.$$
Since
$$\sum_{k=0}^n\f{\bi{2k}{k}^2}{(k+1)16^k}=\f{(2n+1)^2}{16^n(n+1)}\bi{2n}n^2\eq0\pmod{p^2},$$
as observed by van Hamme [vH] (see also (2.5) with $m=16$), we have
$$\sum_{k=0}^n\f{\bi{2k}{k}^2}{(k+1)(-16)^k}\eq-2\sum_{k=0}^n\f{\bi{2k}{k}^2}{(k+1)8^k}\pmod{p^2}.\tag2.7$$

With the help of (2.1),
$$\sum_{k=0}^n\f{\bi{2k}k^2}{(-16)^k}=\sum_{k=0}^n(-1)^k\bi{-1/2}k^2\eq\sum_{k=0}^n(-1)^k\bi nk^2=0\pmod{p}.$$
(Note that $(-1)^{n-k}=-(-1)^k$.) Thus
$$\align&\sum_{h=0}^{p-1}\f{2h+1}{(-16)^h}\sum_{k=0}^h\bi{2k}k^2\bi{2(h-k)}{h-k}^2
\\\eq&\sum_{k=0}^n\f{\bi{2k}k^2}{(-16)^k}\sum_{j=0}^n\f{(2(k+j)+1)\bi{2j}j^2}{(-16)^j}
\eq4\sum_{k=0}^{n}\f{\bi{2k}k^2}{(-16)^k}\sum_{j=0}^n\f{j\bi{2j}j^2}{(-16)^j}\pmod{p^2}.
\endalign$$
By [S5, Lemma 3.2],
$$\sum_{h=0}^{p-1}\f{2h+1}{(-16)^h}\sum_{k=0}^h\bi{2k}k^2\bi{2(h-k)}{h-k}^2\eq p\l(\f{-1}p\r)=-p\pmod{p^2}.$$
Therefore
$$\f1p\sum_{k=0}^n\f{\bi{2k}k^2}{(-16)^k}\sum_{k=0}^n\f{k\bi{2k}k^2}{(-16)^k}\eq-\f14\pmod{p}.\tag2.8$$

In view of (2.5)-(2.8), both (1.5) and (1.6) hold if
$$\sum_{k=0}^n\f{\bi{2k}{k}^2}{(k+1)8^k}\eq\f{(-1)^{(p+1)/4}}2\bi{(p+1)/2}{(p+1)/4}\pmod p.\tag2.9$$

 For $d\in\{0,\}1$,  clearly
$$ a_p^{(d)}(2)=\sum_{x=1}^px^d\l(\f{x(x-1)(x-2)}p\r)=\sum_{r=0}^{p-1}(r+1)^d\l(\f{r(r^2-1)}p\r)$$
and
$$\align a_p^{(d)}(-1)=&\sum_{r=0}^{p-1}r^d\l(\f{r(r^2-1)}p\r)
\\\eq&\sum_{r=0}^{p-1}r^{d+n}(r^2-1)^n=\sum_{k=0}^n\bi nk(-1)^{n-k}\sum_{r=1}^{p-1}r^{n+d+2k}
\\\eq&-\sum\Sb 0\ls k\ls n\\p-1\mid n+d+2k\endSb\bi nk(-1)^{n-k}
\\\eq&\cases 0\pmod{p}&\t{if}\ d=0,
\\(-1)^{(p-3)/4}\bi{n}{(n-1)/2}-\da_{p,3}\pmod{p}&\t{if}\ d=1.\endcases
\endalign$$
Thus we have
$$a_p^{(0)}(2)=a_p^{(0)}(-1)\eq0\pmod{p}$$
and
$$a_p^{(1)}(2)=a_p^{(0)}(-1)+a_p^{(1)}(-1)\eq
(-1)^{(p-3)/4}\bi{n}{(n-1)/2}-\da_{p,3}\pmod{p}.$$
Applying Theorem 1.1 with $\lambda=2$ and $d=0,1$, and noting that
$$\f12\bi{2k+2}{k+1}=\bi{2k+1}{k+1}=2\bi{2k}k-\f{\bi{2k}k}{k+1}\ \ \t{for all}\ k\in\N,$$
we get
$$\sum_{k=0}^n\f{\bi{2k}{k}^2}{(k+1)8^k}+\da_{p,3}\eq 2a_p^{(0)}(2)-a_p^{(1)}(2)\pmod{p}.$$
So (2.9) follows.

(ii) Let $p=2n+1$ be an odd prime. Now we prove (1.7) for all $d\in\{0,\ldots,n\}$ with $d\eq n+1\ (\mo\ 2)$.
(1.7) is valid for $d=n-1$ since
$$\sum_{k=0}^n\f{\bi{2k}k\bi{2k}{k+n-1}}{8^k}=\f{\bi{2(n-1)}{n-1}}{8^{n-1}}+\f{2n\bi{2n}n}{8^n}
=\f{2n+1}{2\times 8^{n-1}}\bi{2n-2}{n-1}\eq0\pmod{p}.$$

Define
$$f(d):=\sum_{k=0}^n\bi{n+k}{2k}\bi{2k}{k+d}(-2)^k\quad \t{for}\ d=0,1,\ldots.$$
Since
$$\bi{n+k}{2k}\eq\f{\bi{2k}k}{(-16)^k}\pmod{p^2}\quad\t{for}\ k=0,\ldots,n\tag2.10$$
(see, e.g., [Su, Lemma 2.2]),
we have
$$f(d)\eq\sum_{k=0}^n\f{\bi{2k}k\bi{2k}{k+d}}{8^k}\pmod{p^2}$$
for all $d=0\ldots,n$.
By applying the Zeilberger algorithm (cf. [PWZ, pp.\,101--119]) via {\tt Mathematica} (version 7), we find the recurrence relation
$$\align &(n-d-1)(n+d+2)(2d+1)f(d+2)
\\=&(2n+1)^2(d+1)f(d+1)-(n-d)(n+d+1)(2d+3)f(d).
\endalign$$
Note that $2n+1=p$. So, if $0\ls d\ls n-2$, then
$$f(d)\eq -\f{(n-d-1)(n+d+2)(2d+1)}{(n-d)(n+d+1)(2d+3)}f(d+2)\pmod{p^2}$$
and hence
$$f(d+2)\eq0\pmod p\ \ \Longrightarrow\ \ f(d)\eq0\pmod p.$$
Now it is clear that (1.7) holds for all $d\in\{0,\ldots,n\}$ with $d\eq n+1\ (\mo\ 2)$.

\heading{3. Proof of Theorem 1.4}\endheading

\medskip
\noindent
{\it Proof of (1.11)}.
Observe that
$$\align \sum_{k=0}^{p-1}\f{\bi{3k}k\bi{2k}k}{27^k}a_k^*
=&\sum_{k=0}^{p-1}\f{\bi{3k}k\bi{2k}k}{27^k}\sum_{m=0}^k\bi km(-1)^ma_m
\\=&\sum_{m=0}^{p-1}(-1)^ma_m\sum_{k=m}^{p-1}\f{\bi{3k}k\bi{2k}k}{27^k}\bi km.
\endalign$$
So it suffices to show that
$$\sum_{k=m}^{p-1}\f{\bi{3k}k\bi{2k}k}{27^k}\bi km\eq\l(\f p3\r)\f{\bi{3m}m\bi{2m}m}{(-27)^m}\ (\mo\ p^2)$$
for all $m=0,1,\ldots,p-1$.

For $0\ls m<n$ define
$$f_n(m)=\sum_{k=m}^{n-1}\f{\bi{3k}k\bi{2k}k}{27^k}\bi km.$$
By Zeilberger's algorithm via {\tt Mathematica} (version 7), we find that
$$\align &9(m+1)^2f_n(m+1)+(3m+1)(3m+2)f_n(m)
\\=&\f{(3n-1)(3n-2)}{27^{n-1}}\bi{n-1}m\bi{2n-2}{n-1}\bi{3n-3}{n-1}.
\endalign$$
Applying this with $n=p>m+1\gs1$ and noting that
$$\bi{2p-2}{p-1}=\f p{2p-1}\prod_{k=1}^{p-1}\f{p+k}k\eq-p\pmod{p^2}\tag3.1$$
and
$$\bi{3p-3}{p-1}=\f{p}{3p-2}\prod_{k=1}^{2p-2}\f{p+k}k\eq-p\pmod{p^2},\tag3.2$$
we get
$$\align &9(m+1)^2f_p(m+1)+(3m+1)(3m+2)f_p(m)
\\\eq&\f{(3p-1)(3p-2)}{27^{p-1}}\bi{p-1}mp^2\eq(-1)^m2p^2\ (\mo\ p^3)
\endalign$$
and hence
$$\align &f_p(m+1)-\l(\f p3\r)\f{\bi{3m+3}{m+1}\bi{2m+2}{m+1}}{(-27)^{m+1}}
\\&+\f{(3m+1)(3m+2)}{9(m+1)^2}\l(f_p(m)-\(\f p3\r)\f{\bi{3m}m\bi{2m}m}{(-27)^m}\)
\\=&f_p(m+1)+\f{(3m+1)(3m+2)}{9(m+1)^2}f_p(m)\eq p^2\f{2(-1)^m}{9(m+1)^2}\ (\mo\ p^3).
\endalign$$
Thus, for every $m=0,\ldots,p-2$, we have
$$\aligned &f_p(m)\eq \l(\f p3\r)\f{\bi{3m}m\bi{2m}m}{(-27)^m}\ \ (\mo\ p^2)
\\ \Longrightarrow\ & f_p(m+1)\eq \l(\f p3\r)\f{\bi{3(m+1)}{m+1}\bi{2(m+1)}{m+1}}{(-27)^{m+1}}\ \ (\mo\ p^2).
\endaligned\tag3.3$$

Since
$$f_p(0)=\sum_{k=0}^{p-1}\f{\bi{3k}k\bi{2k}k}{27^k}\eq\l(\f p3\r)\f{\bi{3\times0}0\bi{2\times0}0}{(-27)^0}\pmod{p^2}$$
by (1.8), with the help of (3.3) we obtain that
$$f_p(m) \eq \l(\f p3\r)\f{\bi{3m}m\bi{2m}m}{(-27)^m}\ (\mo\ p^2)\ \ \ \t{for all}\ m=0,1,\ldots,p-1.$$
This concludes the proof. \qed

\medskip\noindent
{\it Proof of (1.12)}. Similar to the proof of (1.11), we only need to show that
$$\sum_{k=m}^{p-1}\f{\bi{4k}{2k}\bi{2k}k}{64^k}\bi km\eq\l(\f{-2}p\r)\f{\bi{4m}{2m}\bi{2m}m}{(-64)^m}\pmod{p^2}$$
for all $m=0,1,\ldots,p-1$. Since this congruence holds for $m=0$ by (1.9), it suffices to prove that
for any fixed $0\ls m<p-1$ we have
$$\aligned &g_p(m)\eq \l(\f {-2}p\r)\f{\bi{4m}{2m}\bi{2m}m}{(-64)^m}\ (\mo\ p^2)
\\ \Longrightarrow\ & g_p(m+1)\eq \l(\f {-2}p\r)\f{\bi{4(m+1)}{2(m+1)}\bi{2(m+1)}{m+1}}{(-64)^{m+1}}\ (\mo\ p^2),
\endaligned\tag3.4$$
where
$$g_n(m):=\sum_{k=m}^{n-1}\f{\bi{4k}{2k}\bi{2k}k}{64^k}\bi km$$
with $n>m$. By the Zeilberger algorithm, we find that
$$\aligned &16(m+1)^2g_n(m+1)+(4m+1)(4m+3)g_n(m)
\\=&\f{(4n-1)(4n-3)}{64^{n-1}}\bi{n-1}m\bi{2n-2}{n-1}\bi{4n-4}{2n-2}.
\endaligned\tag3.5$$
Clearly
$$\bi{4p-2}{2p-2}=\prod_{k=1}^{2p-2}\f{2p+k}k=\f{3p}p\prod^{2p-2}\Sb k=1\\k\not=p\endSb\l(1+\f{2p}k\r)\eq3\pmod{p}$$
and hence
$$\bi{4p-4}{2p-2}=\f{2p(2p-1)}{(4p-2)(4p-3)}\bi{4p-2}{2p-2}\eq -p\pmod{p^2}.$$
In view of this and (3.1),  from (3.5) with $n=p$ we get
$$16(m+1)^2g_p(m+1)+(4m+1)(4m+3)g_p(m)\eq 3(-1)^mp^2\ (\mo\ p^3).$$
This implies (3.4) since
$$-\f{(4m+1)(4m+3)}{16(m+1)^2}\cdot\f{\bi{4m}{2m}\bi{2m}m}{(-64)^m}=\f{\bi{4(m+1)}{2(m+1)}\bi{2(m+1)}{m+1}}{(-64)^{m+1}}.$$
We are done. \qed

\medskip\noindent
{\it Proof of (1.13)}. For $0\ls m<n$ define
$$h_n(m):=\sum_{k=m}^{n-1}\f{\bi{6k}{3k}\bi{3k}k}{432^k}\bi km.$$
By the Zeilberger algorithm, for $m,n\in\N$ with $m<n-1$, we have
$$\align &36(m+1)^2h_n(m+1)+(6m+1)(6m+5)h_n(m)
\\=&\f{(6n-1)(6n-5)}{432^{n-1}}\bi{n-1}m\bi{3n-3}{n-1}\bi{6n-6}{3n-3}.
\endalign$$
Recall the congruence (3.2) and note that if $p>5$ then
$$\align\bi{6p-6}{3p-3}=&\f{3p(3p-1)(3p-2)}{(6p-3)(6p-4)(6p-5)}\bi{6p-3}{3p-3}
\\\eq&-\f p{10}\prod_{k=1}^{3p-3}\f{3p+k}k\eq-\f p{10}\cdot\f{3p+p}p\cdot\f{3p+2p}{2p}=-p\pmod{p^2}.
\endalign$$
So, no matter $p=5$ or not, for every $m=0,\ldots,p-2$ we have
$$36(m+1)^2h_p(m+1)+(6m+1)(6m+5)h_p(m)\eq0\pmod{p^2}.\tag3.6$$

For $0\ls m<p-1$, since
$$-\f{(6m+1)(6m+5)}{36(m+1)^2}\cdot\f{\bi{6m}{3m}\bi{3m}m}{(-432)^m}=\f{\bi{6(m+1)}{3(m+1)}\bi{3(m+1)}{m+1}}{(-432)^{m+1}},$$
by (3.6) we have
$$\aligned &h_p(m)\eq \l(\f {-1}p\r)\f{\bi{6m}{3m}\bi{3m}m}{(-432)^m}\ (\mo\ p^2)
\\ \Longrightarrow\ & h_p(m+1)\eq \l(\f {-1}p\r)\f{\bi{6(m+1)}{3(m+1)}\bi{3(m+1)}{m+1}}{(-432)^{m+1}}\ (\mo\ p^2).
\endaligned\tag3.7$$
This, together with (1.10), yields that
$$h_p(m)\eq \l(\f {-1}p\r)\f{\bi{6m}{3m}\bi{3m}m}{(-432)^m}\pmod{p^2}$$
for all $m=0,\ldots,p-1$. It follows that
$$\align &\sum_{k=0}^{p-1}\f{\bi{6k}{3k}\bi{3k}k}{432^k}\sum_{m=0}^k\bi km(-1)^ma_m
\\=&\sum_{m=0}^{p-1}(-1)^ma_mh_p(m)\eq \l(\f{-1}p\r)\sum_{m=0}^{p-1}a_m\f{\bi{6m}{3m}\bi{3m}m}{432^m}\ (\mo\ p^2).
\endalign$$
This proves (1.13). \qed

\heading{4. Proofs of Theorems 1.5--1.6}\endheading

\medskip
\noindent{\it Proof of Theorem 1.5}. By (1.23),
$$\sum_{k=0}^{p-1}\f{\bi{4k}{2k}\bi{2k}k}{(k+1)72^k}\eq p+\f{72-64}{12}\sum_{k=0}^{p-1}\f{k\bi{4k}{2k}\bi{2k}k}{72^k}\pmod{p^2}$$
and hence
$$\sum_{k=0}^{p-1}\f{k\bi{4k}{2k}\bi{2k}k}{72^k}\eq\f 32\sum_{k=0}^{p-1}\f{\bi{4k}{2k}\bi{2k}k}{(k+1)72^k}\pmod{p}.$$
So it suffices to determine $\sum_{k=0}^{n}k\bi{4k}{2k}\bi{2k}k/{72^k}$ mod $p$, where $n=(p-1)/2$.
(Note that $p\mid\bi{2k}k$ for $k=n+1,\ldots,p-1$.)

The Legendre polynomial of degree $n$ is given by
$$P_n(x):=\sum_{k=0}^n\bi nk\bi{n+k}k\l(\f{x-1}2\r)^k=\sum_{k=0}^n\bi {n+k}{2k}\bi{2k}k\l(\f{x-1}2\r)^k.$$
It is known (see, e.g., [N]) that
$$\sum_{k=0}^n\bi n{2k}\bi{2k}kx^k=(\sqrt{1-4x})^nP_n\l(\f1{\sqrt{1-4x}}\r).$$
Taking derivatives of both sides of this identity, we get
$$\align&\sum_{k=0}^n\bi n{2k}\bi{2k}kkx^{k-1}
\\=&-2n(1-4x)^{n/2-1}\sum_{k=0}^n\bi {n+k}{2k}\bi{2k}k\l(\f{(1-4x)^{-1/2}-1}2\r)^k
\\&+(1-4x)^{(n-3)/2}\sum_{k=0}^n\bi{n+k}{2k}\bi{2k}kk\l(\f{(1-4x)^{-1/2}-1}2\r)^{k-1}.
\endalign$$
In view of (2.1) and (2.10), by putting $x=2/9$ in the last equality we obtain
$$\f12\sum_{k=0}^n\f{k\bi{4k}{2k}\bi{2k}k}{72^k}
\eq\f 1{3^n}\sum_{k=0}^n\f{\bi{2k}k^2}{(-16)^k}+\f{3}{3^n}\sum_{k=0}^n\f{k\bi{2k}k^2}{(-16)^k}\pmod p$$
and hence
$$\l(\f{3}p\r)\sum_{k=0}^n\f{k\bi{4k}{2k}\bi{2k}k}{72^k}
\eq2\sum_{k=0}^n\f{\bi{2k}k^2}{(-16)^k}+6\sum_{k=0}^n\f{k\bi{2k}k^2}{(-16)^k}\pmod{p}.$$
Since $\sum_{k=0}^n\bi{2k}k^2/(-16)^k$ and  $\sum_{k=0}^nk\bi{2k}k^2/(-16)^k$ modulo $p$ have been determined
(cf. Theorem 1.3(i) and the paragraph after (1.4)),  we finally obtain the desired result for $\sum_{k=0}^nk\bi{4k}{2k}\bi{2k}k/72^k$ modulo $p$.

The proof of Theorem 1.5 is now complete. \qed

\medskip\noindent{\it Proof of Theorem 1.6}.
We just show the first part in detail. Parts (ii) and (iii) can be proved similarly.

By (1.14) and (1.17), there are $p$-adic integers $a_0,\ldots,a_{p-1},b_0,\ldots,b_{p-2}$ such that
$$\sum_{k=0}^{p-1}\f{\bi{3k}k\bi{2k}k}{27^k}\l(x^k-\l(\f p3\r)(1-x)^k\r)=p^2\sum_{k=0}^{p-1}a_kx^k\tag4.1$$
and
$$\sum_{k=1}^{p-1}\f{k\bi{3k}k\bi{2k}k}{27^k}\l(x^{k-1}+\l(\f p3\r)(1-x)^{k-1}\r)=p^2\sum_{k=0}^{p-2}b_kx^k.\tag4.2$$
Let $\al$ and $\beta$ be the two distinct roots of the equation $x^2-Ax+B=0$. It is well known that
$$u_k=\f{\al^k-\beta^k}{\al-\beta}\ \ \t{and}\ \ v_k=\al^k+\beta^k\quad \ \t{for all}\ k\in\N.$$
As $\al/A+\beta/A=1$, by (4.1) and (4.2) we have
$$\align\sum_{k=0}^{p-1}\f{\bi{3k}k\bi{2k}k}{(27A)^k}\l(u_k+\l(\f p3\r)u_k\r)&=p^2\sum_{k=0}^{p-1}\f{a_k}{A^k}u_k,
\\\sum_{k=0}^{p-1}\f{\bi{3k}k\bi{2k}k}{(27A)^k}\l(v_k-\l(\f p3\r)v_k\r)&=p^2\sum_{k=0}^{p-1}\f{a_k}{A^k}v_k,
\\\sum_{k=1}^{p-1}\f{k\bi{3k}k\bi{2k}k}{27^kA^{k-1}}\l(u_{k-1}-\l(\f p3\r)u_{k-1}\r)&=p^2\sum_{k=0}^{p-2}\f{b_k}{A^k}u_k,
\\\sum_{k=1}^{p-1}\f{k\bi{3k}k\bi{2k}k}{27^kA^{k-1}}\l(v_{k-1}+\l(\f p3\r)v_{k-1}\r)&=p^2\sum_{k=0}^{p-2}\f{b_k}{A^k}v_k.
\endalign$$
Thus (1.30) holds when $p\eq1\ (\mo\ 3)$, and (1.31) holds when $p\eq2\ (\mo\ 3)$. We are done. \qed

\medskip

\Ack. The author is grateful to the referee for helpful comments.

\enddocument